\documentclass[12pt]{article}

\usepackage{t1enc}
\usepackage[latin1]{inputenc}
\usepackage[french]{babel}

\usepackage{amsmath}
\usepackage{amssymb}

\usepackage{colortbl}

\usepackage[T1]{fontenc} 
\usepackage{times}
\usepackage{graphicx}

\newtheorem{Def}{Definition}[section]
\newtheorem{Pro}{Proposition}[section]
\newtheorem{Cor}{Corollary}[section]
\newtheorem{Rem}{Remark}[section]
\newtheorem{Lem}{Lemma}[section]
\newcommand{\ds}{\displaystyle}

\title{A tree-growth model to optimize silviculture}

\author{Patrice Loisel\thanks{ INRA, UMR 729 MISTEA, 2 place Viala, F-34060 Montpellier, France}
  \thanks{ SupAgro, UMR 729 MISTEA, 2 place Viala, F-34060 Montpellier, France}, Jean François Dhôte\thanks{Office National des Forêts, R\&D Departement
Boulevard de Constance, F-77300 Fontainebleau, France}}

\date{2011}

\begin{document}
\baselineskip=16pt
\maketitle{}

{\bf Abstract}: In this paper, we present the description of a
simplified model of the dynamic of a mono-specific even-aged
forest. The model studied is a tree-growth model based on a system of
two ordinary differential equations concerning the tree basal area and
the number of trees. The analytical study of this model permits us to
predict the behavior of the system solutions. We are trying to
highlight the influence of economic parameters and growth parameters
on the system solutions, in the framework of the optimization of
silviculture.

\section{Introduction}
The forest management, because of its impact on our environment, is a
topic that now involves researchers of many disciplines: forestry,
economy, and ecology. These various communities have models adapted to
the questions they wish to tackle. 
The economists usually study the best age at which to cut down a tree
or stand of trees, the simultaneous management of several forest
stands.
 As to foresters, they are moreover interested in
silviculture at the stand level. We will focus on the models developed
by foresters.

The models of growth for silviculture, expanded rapidly these last few
years, and represent a significant part of the developed models.
Here, we are focusing on a particular type of forest (thus a
particular type of model): a mono-specific even-aged forest where all
trees belong to the same species and are the same age.

Models built by forest modellers are based on statistical adjustments
of dendrometric data \cite{Porte}: these models accurately describe
the evolution of a forest, but the analytical study of those models is
made difficult due to their complexity. The analytical study of models
allows us to predict the influence of different parameters: for
instance, if by modifying parameters to take into account the climate
change and analyzing so the potential consequences is provide. To
allow analytical studies while having realistic model, we decide in
this paper to consider a simplified model.

 The models, on which we are focusing here, are tree centered distance
independent models where the trees are not spatialized.  In this type
of model each tree is characterized by its basal area at the height of
$1.3$ meters: $s$ and eventually by its height $h$. The model here
described is based on the concepts developed in the growth model
``Fagacées'' \cite{Jef} \cite{Gilles} for Oak or Beech forest. In this
model, the link between the stand level and the individual tree level
is explicit. This modeling allows us to describe the evolution of a
forest of high density.  ``Fagacées'' was broadcast through the
project Capsis \cite{Capsis}.

In order to allow analytical studies, we are starting with the
simplified hypothesis in which we consider that all trees have the
same basal area. We then consider a forest of $n$ trees with basal area $s$
which needs management such as thinning $e$ throughout time.

The tree growth (due to the observed densities) is not independent of its
neighbor's growth.  There is a competition for the available
resources: photosynthesis and access to the light on one hand, and
mineral nutrients on the other hand. Thus, cutting a tree implies the
increase of its neighbors' growth, and cutting no tree limits
individual growth.  That shows how a forest is not just a juxtaposition
of trees.

This phenomenon is considered through an assessment equation that
allows us to distribute the energy resources of the forest stand level
between the various trees. This equation coupled to an equation
describing the evolution of the number of trees leads us to a dynamic
system of the forest.  The studied model takes into account the
characteristics that the foresters consider as the most important and
as required: the basal area at $1.3$ meters, the height, and the number
of trees.

In Section 2 we will present the designing of the model. In Section 3
we will present general results on the behavior of the system
solutions, then we will look for strategies which permit to leave the
viability domain in minimal or maximal time. Finally in Section 4 we
will highlight the influence of economic parameters and growth
parameters for silviculture within the resolution of an optimization
problem.

\section{Designing the model}
{\bf The trees density and the RDI of forest stand}.
 
Let's consider a forest with a given area. It is intuitively clear
that the tree number which this area can bear is limited. The
environmental conditions (type of soil, local climatic conditions) are
also factors to be taken into account for the maximum tree capacity.  Foresters
have established a law called ``self-thinning'', described hereafter,
to evaluate this maximum capacity. Let's note $s_*$ the average
tree basal area (at the height of $1.3$ meters) of the forest.      Reineke
\cite{Reineke} observed monospecific forests with various densities and
various species.   Out of these observations he claimed the maximum
tree number $n_{max}(s_*)$ that a stand can bear is given by the
following self thinning relation:

 $$\log n_{max}(s_*) = C_0 - {q \over 2} \log s_*,$$

where $C_0 > 0$ and $1 < q < 2$ are characteristic constant
values of the forest species and of its environment, and in particular the
ground fertility.

As for a given $s_*$, beyond the number $n_{max}(s_*)$ the trees die,
the forest stand (in terms of tree number) has to remain under this
limit. To simplify we will here make the assumption that all the trees
have the same basal area $s$.  For a forest of where the effective tree
number is $n$, taking into account the relation of self-thinning, the
density $r$ is the ratio of the tree number and the maximum tree
number of basal area s that the forest can support, $r$ is then defined
by:

 $$r(n,s)={n \over   n_{max}(s)}$$

 is written this way:

$$r(n,s) = {n s^{q \over 2} \over e^{C_0}} = A n s^{q \over 2}$$
 where  $A:=e^{-C_0}$. This ratio $r$ is called $RDI$ (Relative
 Density Index or Reineke Density Index). By definition this ratio is
 always less than $1$.   

{\bf Competition between trees: from forest stand level to individual
tree level}

We will now describe the temporal evolution of state variables $s, n$
and now $r$, of the considered forest.

The growth of a tree depends on its neighbors. There is competition
for the resources and the death of a tree, natural or due to cuttings,
implies an increased growth for its neighbors .  We make the
assumption, that in the course of time, silviculture makes it possible
to maintain the trees uniformly distributed on the area. The model 
is characterized by the existence of two levels in the modeling.

{\it At the forest stand level} the available energy for the
considered forest, is considered globally, due to photosynthesis or
due to nutrients in the soil.  This supplied energy makes it possible
to ensure at the same time the maintenance and the growth of the
trees. The share reserved for maintenance increases with the tree
height, therefore with time, which limits all the more so the
available part for growth.  The energy left for growth is therefore a
decreasing time function and allows the increase of basal area of the
forest stand.  The increase of basal area of the forest stand at its
peak of density $r(n(t),s(t))=1$ is given by the function $V (.)$. We
assume that $V (t)$ verify the following properties:

$(H_1)$: $V(.)$ is a positive, decreasing, convex function of $t$.

For a lower density ($r(n(t),s(t)) < 1$), the effective increase of
the basal area for the forest is reduced by a factor dependent on this
same density: $g(r(n(t),s(t)))$ at any time.  Thus the energy actually
used at time $t$ is given by:

$$g(r(n(t),s(t))) V(t).$$

The function $g(.)$ is supposed to satisfy the following properties:

$(H_2)$: $g(.)$ is an increasing, concave function of $r$ such that
$g(r) > r$ for  $r \in (0,1)$, $g(0)=0, g(1)=1$.

The concavity of $g$ is related to crown development in relation to
basal area. 

{\it On the individual tree level}, tree growth is characterized by
the evolution of tree basal area and therefore by the evolution of the
function $s(t)$: the instantaneous increase is thus $\ds {ds(t) \over
dt}$.  As mentioned in the hypothesis all trees have the same basal area,
the total sum increase of basal areas of all trees is $\ds n(t) {ds(t)
\over dt}$.  This total increase is obtained from the available energy
resources.  We thus obtain the equation which describes the link
between the forest stand level and the individual tree level:

$$ g(r(n(t),s(t))) V(t) = n(t) {ds(t) \over dt } \mbox{ for } n(t) > 0, \mbox{
for all } t $$

For any $t$, this enables us to establish the first dynamic equation
of our model:
$${ds(t) \over dt } = {g(r(n(t),s(t))) \over n(t)} V(t) $$ In addition, the
evolution of the tree number depends on several factors. To permit
analytical study of the model, we've decided to simplify and we
suppose the only cause of tree mortality is due to fallings that
foresters could operate.   We noted $e(t)$ the instantaneous rate of
trees cutting at time $t$.  Thus the evolution of the tree number
is given by:

$$\ds {dn(t) \over dt} = - e(t).$$

We wish to preserve a minimum tree number in the forest stand, which
implies $n(t) \geq \underline n > 0$, for any $t$. Technologically and
to ensure a provisioning not too irregular, the thinning rate is
limited: $0 \leq e(t) \leq \overline e$, for any $t$.

The forest is therefore described using the two state variables $s, n$
and its evolution follows the following dynamic:

\begin{equation*}
({\mathcal S_0}) \ \ \ 
\begin{cases}
\ds {ds(t) \over dt} = {g(r(n(t),s(t))) \over n(t)} V(t) \\
\ds {dn(t) \over dt} = -e(t) 
\end{cases}
\end{equation*}

with the constraints $0 \leq e(t) \leq \overline e$, $n(t) \geq
\underline n, r(n(t),s(t))=A n(t) s(t)^{q \over 2} \leq 1$ for any
$t$.

Foresters built this type of model from observed forest data. The 
available data only allows us to validate the model on a limited period of
time.  The system $({\mathcal S_0})$ has therefore a time limit domain:
$t \in [0,T_*]$.

It is a dynamic system in the state variables $n$ and $s$, controlled
by the control variable $e$. For a cutting policy, i.e. the data of a
particular function $e(.)$, and for each initial condition $(s(0),
n(0))$, this system has a single solution: we will suppose the
functions $g(.), V (.)$ are regular enough for it to happen. We will
specify these trajectories in the following paragraph.

To finish with the model description, the tree height $h$ is supposed
to depend only on the tree basal area $s$ and on the dominant height $h_0$
(average height of the 100 largest trees), $h_0$ is a concave function
of time $t$ and shouldn't depend on silviculture (cuttings in the
course of time) and thus depends only on time $t$. The height $h$ has
therefore no influence on $s(t)$ and its evolution, $h$ is
consequently an output of the model.

If, as supposed earlier, the basal area $s$ at time $t$ is the same for all trees,
the height $h$ is also the same. We therefore deduce $h(t) = h_0(t)$, for any $t$.

\section{Studiing the solutions}

\subsection{Model properties}

The solutions of the dynamic system $({\cal S}_0)$ must satisfy in
particular the constraint $r(n(t),s(t)) \leq 1$ for any $t$. If there
is no cutting, i.e. if $e (t) = 0$ for any $t$, we deduce that $n (t)
= n (0)$, $s (.)$ and $r (.)$ are increasing with respect to time $t$.
   Let's suppose there is one time $\tau < T_*$ such as
$r(n(\tau),s(\tau))=1$, we deduced $t > \tau$  if we apply a control
identically null then $r(n(t),s(t)) > 1$ and the constraint is no
longer satisfied. In order to let us know which control we should
apply we are led to study the evolution of the density function
$r(.)$:

\begin{align}
  \label{r}
  {dr(n(t),s(t)) \over dt} & \ds = r'_s(n(t),s(t)) {g(r(n(t),s(t))) \over n(t)}
 V(t) -r'_n(n(t),s(t)) e(t) \notag \\
& \ds =  {r(n(t),s(t)) \over n(t)} [{q \over 2}{g(r(n(t),s(t))) \over s(t)}V(t) -e(t)]
\end{align}

Out of this last equation we can deduce that in order to respect the
constraint $r(n(t),s(t)) \leq 1$ for $t > \tau$, we should apply a
non-identically null control on the system.  Thus the cutting $\ds
e(t)={q \over 2}{V(t)\over s(t)}$ for $t > \tau$ respects the
constraints by binding, i.e. $r(n(t),s(t)) = 1$. If we define the
function $e_r(., .)$ by: $\ds e_r(s,t):= {q \over 2}{V(t)\over s}$,
for any $s > 0$, $t > 0$, the solutions, independently of the cutting
function $e(.)$ applied to the system, are only valid if the
constraint: $\ds e_r(s(t),t) \leq \overline e$ is satisfied. We are
therefore led to formulate the following assumption $(H_3) $:

$(H_3)$: $\ds e_r(s_m(t),t)={q \over 2} {V(t) \over s_m(t)} <
\overline e$ for all $t \in (0,T_*)$

where $s_m (t)$ is the minimal value $s (t)$ can reach at the time $t$.

\begin{Rem}
\label{sm} 
$s_m(t)$ is not specified at the moment but will be specified later
on, however we can take an approximate lower bound for now: $s_m(t) >
s(0)$.
\end{Rem}

We noted previously that the system $(\mathcal S_0)$ is considered
only for $t \in [0,T_*]$. It is advisable to specify now, the behavior
of the solutions in this interval.

\begin{Def}  The function $\cal V(.;.)$ is defined by
$\ds {\mathcal V}(t;T)=\int_t^{T} V(u) du$ and represents the energy
that has been available for growth in the period  $[t,T]$.
\end{Def}

The following Lemma shows us that the system validity field depends on
this energy value:

\begin{Lem}
\label{lemval} Assuming $(H_2), (H_3)$, then:

(i) if $\ds {\mathcal V}(0;T_*) $ is large enough then there exists a
time $\tau < T_*$ such that $r(n(\tau),s(\tau))=1$ and $n(\tau)=\underline
n$. The dynamical system is only valid on the interval $[0,\tau]$.
This time $\tau$ depends on the evolution of the cutting $e(.)$.

(ii) conversely if $\ds {\mathcal V}(0;T_*) $ is small enough then the
dynamical system is valid throughout the entire interval $[0,T_*]$.

\end{Lem}

{\it Proof}: (i) From $g(r) \geq r$ we deduce: $\ds {ds(t) \over dt}
\geq {r(n(t),s(t)) \over n(t)} V(t) = A s(t)^{q \over 2} V(t)$ hence:
$$ s(T_*)^{1-{q \over 2}} \geq s(0)^{1-{q \over 2}} + A(1-{q \over 2})
{\mathcal V}(0;T_*) $$

From $r(n(T_*),s(T_*)) \leq 1$ and $n(T_*) \geq \underline n$ we
deduce: $\ds s(T_*) \leq {1 \over (A \underline n)^{2 \over q}}$. If
$\ds {\mathcal V}(0;T_*) $ is large enough, we obtain a contradiction.

  (ii) Let's set $\tau$ the first period where $r(n(t),s(t))$ reaches
  $1$, then for $0 < t \leq \tau $ we deduce: $\ds {ds(t) \over dt}
  \leq {V(t) \over n(t)} \leq {V(t) \over \underline n}$ and therefore
  $\ds s(\tau) \leq s(0) + { {\mathcal V}(0;\tau) \over \underline
    n}$.  If $\ds {\mathcal V}(0;T_*) $ is small enough, we deduce
  $r(n(\tau),s(\tau)) < 1$ in contradiction with the
  assumption. \hfill $\square$

    Specific trajectories easily expressed in terms of control, will
  play an important role, we are introducing them here: let's consider
  the system of equations $({\mathcal S_0})$, for trajectory $E_0$
  from a fixed initial condition $(s(0), n(0))$ we apply the maximum
  cutting $e(t) = \overline e$ until we reach the value $n$ for the
  tree number, $t_{0,n}$ is the time needed to go from the tree number
  $n(0)$ to $n$. By definition, we therefore have $\ds t_{0,n} =
  {n(0)-n \over \overline e}$.

  For trajectory $E^0$, starting from the same initial condition (with
  $r(n(0),s(0)) < 1$) we apply the minimum cutting $e(t) = 0$ until
  reaching the value $1$ for the $RDI$ $r$, then we apply the control
  $e_r(s(t),t)$ until $n = \underline n$. $t^0$ is the time needed to
  go for the $RDI$ from $r(n(0),s(0))$ to $1$ and $T^0$ the final
  time. By definition, we therefore have $t^0$ and $T^0$ respectively
  solutions of:

$$ {q \over 2} n(0)^{{2 \over q}-1} A^{2 \over q} {\mathcal V}(0;t^0) =
\int_{r(n(0),s(0))}^1 {u^{{2 \over q}-1} \over g(u)} du \mbox{ (at
  constant $n$) }$$

$$\underline n^{1-{2 \over q}}=n(0)^{1-{2 \over q}} 
+A^{2 \over q} (1-{q \over 2}) \mathcal V(t^0;T^0) \mbox{ (at constant
  $r$) }$$

{\bf Notations.} To summarize we note the following definitions of the
specific trajectories:

\noindent \begin{align}
E_{0}: e(t) & = 
\begin{cases} 
\overline e & \mbox{ if } t < t_{0,\underline n} \\
 0 \ \ \ \ \ \ \ \ \ \ \ \ \ \  & \mbox{ if }  t > t_{0,\underline n}
\end{cases} \notag \\
E^0: e(t) & =
\begin{cases} 
0 &  \mbox{ if } r(n(t),s(t)) <1, \mbox{ i.e. } t < t^0 \\
\ds e_r(s(t),t) & \mbox{ if } t^0 < t < T^0
\end{cases} \notag \\
\mbox{ For  } t_{0,\underline n} < T < t^0 &
 \mbox{ we can also define an intermediate trajectory } E_T \notag \\
E_T: e(t) &=
\begin{cases} 
0 \ \ \ \ \ \ \ \ \ \ \ \ \ \ \   & \mbox{ if } t < T-t_{0,\underline n} \\
\overline e & \mbox{ if } t > T-t_{0,\underline n}
\end{cases} \notag \\
& \mbox{ if } t_{0,\underline n} < T <t_{0,\underline n}+t^0 \notag \\
E_T: e(t) &=
\begin{cases} 
0 \ \ \ \ \  & \mbox{ if } t < t^0 \\
e_r(s(t),t) & \mbox{ if } t^0 < t < t_* \\
\overline e &  \mbox{ if } t_* < t < T
\end{cases} \notag \\
& \mbox{ if } t_{0,\underline n} +t^0  < T <  T^0 \notag
\end{align}

where $t_*$ is defined by: $\ds n(t_*)^{1-{2 \over q}}=n(0)^{1-{q
    \over 2}} +A(1-{2 \over q}) \mathcal V(t^0;t_*)$ and $(T-t_*)
\overline e=n(t_*) - \underline n$.

We note that $E_{T^0} = E^0$ and by extension if $T < t_{0,\underline
  n}$ then $E_T=E_0$.

\begin{center}

Figure 1: Phase plane in the coordinates $s$ and $n$.
\end{center}

The functions obtained by just following the trajectories $E_{0}, E_T$
and $E^0$ will be noted by the indices $_{0}$, $_T$ and $^0$.

From the increasing of the basal area $s$ and the non increasing of the
number $n$ we deduce that the system solutions have no choice but to
move to the right bottom in the phase plane.  The Lemma \ref{lemval}
(i) has shown that, if $\ds {\mathcal V}(0;T_*) $ is large enough, the
solutions are not valid throughout the entire interval $[0,T_*]$. That
implies that as from a time $\tau$, the solution doesn't belong to the
validity domain defined by the constraints $r(n(t),s(t)) \leq 1$ and
$n(t) \geq \underline n$. The only point which makes it possible to
leave this validity domain is the point such as $r(n(\tau),s(\tau))=1$
and $n(\tau) = \underline n$,  this point is represented by a square
on  Figure 1. As the solution remains valid basal area $s(t)$ verifies
for all $t < T$: $\ds s(t) \leq \overline s = {1 \over (A \underline
  n)^{2 \over q}}$ (is deduced from $r(n(t),s(t)) \leq 1$).

We define $\underline T$ (resp. $\overline T$) as the minimum
(resp. maximum) time necessary to reach the point defined by
$r(n(T),s(T)) = 1$ and $n(T) = \underline n$. Then:

- if $T \leq \underline T$ the solution remains valid whatever the
trajectory (i.e.  whatever the evolution of the cutting $e(.)$).

- if $\underline T < T \leq \overline T$  the system has a solution on $[0,T]$ for certain controls $e(.)$.

- if $T >\overline T$ the system has no solution on $[0,T]$ whatever the controls $e(.)$.

We noted that, from Lemma \ref{lemval}, if $\ds {\mathcal V}(0;T_*)$
is small enough then $\underline T$ and especially $\overline T$ can
no exist.

{\bf A particular case}

We could consider the particular function $\ds g(r) = g_{\theta}(r) =
r^{1-\theta}, 0 < \theta < 1$.  In that case $\ds {ds(t) \over dt} =
A^{1-\theta} {s(t)^{{q \over 2}(1-\theta)} \over
  n(t)^{\theta}}V(t)$. The basal area $s$ is explicitly deduced from the
tree number $n$:

$$ s(t)^{1-{q \over 2}(1-\theta)} = s(0)^{1-{q \over 2}(1-\theta)}
+ A^{1-\theta}(1-{q \over 2}(1-\theta)) \int_0^t {V(u) \over
  n(u)^{\theta}}du.$$

 In this class of functions $g_{\theta}(.)$ we consider the extreme
 case ($\theta=0$) for which some of the properties of the hypothesis
 $(H_2)$ are not satisfied: $g_0(r) = r$. In this last case ${\cal G}
 (r) \equiv 0$ and the evolution of the basal area $s$ is independent
 from the evolution of the tree number $n$:

$$ s(t)^{1- {q \over 2}} =s(0)^{1- {q \over 2}} +A (1- {q \over 2})
{\mathcal V}(0;t).$$

In this particular case, provided that $\ds {\mathcal V}(0;T_*)$ is
large enough, $\underline T$ and $\overline T$ are equal, don't depend
on the cutting $e(.)$ and are the unique solution of the following
equation in $T$:

$$s(0)^{1- {q \over 2}} +A (1- {q \over 2}) {\mathcal V}(0;T)
= \overline s^{1-{q \over 2}}.$$

\subsection{Minimum and maximum time necessary to reach the point
  $(r,n)=(1,\underline n)$}

In order to succeed in the conclusion of the study of the minimum and
maximum time needed to reach this point, we will need the following
properties and definitions related to the function $g(.)$. The
increase in basal area of each tree is given by $\ds {g(r(n(t),s(t)))
  \over n(t)} V(t)$.  We will to know thereafter the evolution of $\ds
{g(r(n(t),s(t))) \over n(t)}$, for the same aim, we will need to
define the functions ${\cal G}(.)$, $\gamma(.)$:

\begin{Def}  
  The function $\cal G(.)$ is defined by: $\ds {\cal G}(r)= {d \over
    dr}[{r \over g(r)}]= {g(r)-rg'(r) \over g^2(r)}$. 
  The function
  $\gamma(.)$ is defined by: $\ds \gamma(r)={r g'(r) \over g(r)}$.
\end{Def}

In the ``Fagacées'' model, $\ds g(r) = {(1+p)r \over r+p}$ with $p >
0$, ${\cal G}$ is constant $\ds {\cal G}(r) \equiv G = {1 \over 1+p} $
.

From ${\cal G}(.)$ and $\gamma(.)$ definitions, we can establish the
following properties for the model:

\begin{Lem} 
\label{lemG}
 Assuming the hypothesis $(H_2)$, then:

 (i) The function $\ds {r \over g(r)}$ is an increasing function of
 $r$ and $\cal G(.)$ satisfies $0 < {\cal G}(r) g(r) \leq 1$ for any
 $r > 0$.

(ii) The function $\ds {g(r(n,s)) \over n}$ is a decreasing function of $n$.

(iii) The function $\ds g(r(n,s))$ is an increasing function of $s$.

(iv) The function $\gamma(.)$ satisfies $\gamma(r) \leq 1$ for any $r
> 0$.

(v) if $n(t)$ and $s(t)$ are solutions of systems $({\mathcal S_0})$
the function $\ds {g(r(n(t),s(t))) \over n(t)}$ is an increasing
function of $t$.
\end{Lem}

{\it Proof}: (i) From the concavity of $g(.)$, $\ds {d \over dr}[g(r)
- rg'(r)] = -r g''(r) > 0$ for any $r > 0$ and from $g(0) = 0$, we
deduce that $\ds g(r) -r g'(r) > 0$ and ${\cal G}(r) > 0$.  From $\ds
g'(r) > 0$, $\ds {\cal G}(r) \leq {1 \over g(r)}$.

(ii) $\ds {\partial \over \partial n}[{g(r(n,s)) \over n}] 
= {r(n,s)g'(r(n,s)) -g(r(n,s)) \over n^2}= - {\cal
  G}(r(n,s)){g^2(r(n,s)) \over n^2} <0$.

(iii) $\ds {\partial g(r(n,s)) \over \partial s} = g'(r(n,s))
r'_s(n,s) > 0$.

(iv) { From $ \ds g(r) -r g'(r) > 0$ we deduce the result}.

(v) $\ds {d \over dt} [{g(r(n(t),s(t))) \over n(t)}] =
{(g'(r)r'_s)(n(t),s(t))) \over n(t)} {ds(t) \over dt} + {\cal
  G}(r(n(t),s(t))){g^2(r(n(t),s(t))) \over n^2(t)} e(t)$

and from (i) we deduce the result. \hfill $\square$

We are thus focusing on the trajectories and also on the strategies
which allow us to reach respectively in a minimum and maximum time the
point $(1,\underline n)$ in the $(r,n)$ coordinates.

The minimum time $\underline T$ (resp. the maximum time $\overline T$)
is reached by solving the problem of optimal control: $\ds \min_{e(.)}
T$ (resp. $\ds \max_{e(.)} T$) with the set of admissible values for
the control variable $[0,\overline e]$. $n(.)$ and $s(.)$ are the
state variables governed by the system $({\mathcal S_0})$ of initial
condition $(n(0), s(0))$ and satisfying constraints, for all $t \in
[0,T)$, $r(n(t),s(t)) \leq 1$, $n(t) \geq \underline n$ and { the
  right end time} constraints $r(n(T),s(T)) =1, n(T)=\underline n$.



\begin{Pro}
\label{mini}
Assume $(H_2), (H_3)$. If $n(.)$ and $s(.)$ are the solutions of the
system $(\cal S_0)$ for a control $e(.)$ then, $ \forall t \in [0,T]$:

(i) $\ds n_0(t) \leq n(t) \leq n^0(t)$ 

(ii) $ s(t) \geq s^0(t) $

(iii) if $g(r) = g_{\theta}(r)=r^{1-\theta}$, $ s(t) \leq s_0(t)$

If we assume that the final tree-number $n(T)$ is equal to $\underline
n$ then:

(iv) $\ds n(t) \leq n_T(t)$ 

(v) if $g(r) = g_{\theta}(r)=r^{1-\theta}$, $ s(t) \leq s_T(t)$

\end{Pro}

{\it Proof} (i) and (iv) Follows from the definition.

(ii) From Lemma \ref{lemG} (ii) we deduce: $\ds {ds(t) \over dt} =
{g(r(n(t),s(t)) \over n(t)} V(t) \geq {g(r(n^0(t),s(t)) \over n^0(t)}
V(t)$.

For $t \leq t^0$ we deduce: $\ds {ds \over g(r(n(0), s^{q\over 2}))}
\geq {V(t) \over n(0)} dt$ then by integration of the inequality:

$\ds \int_{s(0)}^{s(t)} {dx \over g(r(n(0),x^{q\over 2}))}$
$\ds \geq {\mathcal V(0;t) \over n(0)} = \int_{s(0)}^{s^0(t)} {dx
  \over g(r(n(0),x^{q\over 2}))}$
and we deduce $s(t) \geq s^0(t)$.

For $ t > t^0$, $\ds {ds(t) \over dt} \geq {V(t) \over n^0(t)} =
{ds^0(t) \over dt}$ and from $s(t^0) \geq s^0(t^0)$ we deduce by
integration $s(t) \geq s^0(t)$.

(iii) From the previously stated expression of the basal area $s$ and
$n(t) \geq n_0(t)$ we deduce:




\begin{align}
  {s(t)^{1-{q \over 2}(1-\theta)}-s(0)^{1-{q \over 2}(1-\theta)} \over
    A^{1-\theta} (1-{q \over 2}(1-\theta))} = \int_0^t {V(u) \over
    n(u)^{\theta}} du \leq \int_0^t {V(u) \over n_0(u)^{\theta}} du=
  {s_0(t)^{1-{q \over 2}(1-\theta)}-s(0)^{1-{q \over
        2}(1-\theta)}\over A^{1-\theta} (1-{q \over 2}(1-\theta))}
  \notag
\end{align}

and hence the result. 

(v) From the previously stated expression of the basal area $s$ and
$n(t) \leq n_T(t)$ we deduce:




\begin{align}
  {s(t)^{1-{q \over 2}(1-\theta)}-s(0)^{1-{q \over 2}(1-\theta)} \over
    A^{1-\theta} (1-{q \over 2}(1-\theta))} = \int_0^t {V(u) \over
    n(u)^{\theta}} du \geq \int_0^t {V(u) \over n_T(u)^{\theta}} du=
  {s_T(t)^{1-{q \over 2}(1-\theta)}-s(0)^{1-{q \over
        2}(1-\theta)}\over A^{1-\theta} (1-{q \over 2}(1-\theta))}  \notag
\end{align}

\hfill $\square$

 
If we remark that to reach the point $(r,n)=(1,\underline n)$ in
minimal time (resp. in maximal time) is equivalent to reach $s=
\overline s$ in minimal time (resp. in maximal time), we deduce the
the trajectory that allows to reach the point $(r,n)=(1,\underline n)$
in minimal or maximal time:

\begin{Cor}
\label{minimum}
Assume $(H_2), (H_3) $. Let $\underline T$ the minimal time
(resp. $\overline T$ the maximal time) necessary to reach the point
$(r,n)=(1,\underline n)$ using the
control $e(.)$. Then:

(i) if $\underline T$ is finite and $g(r) = r^{1-\theta}$, the
trajectory that allows to reach the point $(r,n)=(1,\underline n)$ in
minimal time $\underline T$ is the trajectory $E_{0}$.

(ii) if $\overline T $ is finite, the trajectory that allows to reach
the point $(r,n)=(1,\underline n)$ in maximal time $\overline T$ is
the trajectory $E^0$.  Maximal time $\overline T$ is then the solution
of:

$$\underline n^{1-{2 \over q}}-n(0)^{1-{2 \over q}}=
A^{2 \over q} (1-{q \over 2}) {\mathcal V}(t^0;\overline T).$$
\end{Cor}



\section{Optimization of silviculture}

{ In order to optimize the silviculture, we are interested in
  problems which consist in seeking the minimal and maximum values of
  a variable function depending on the state variables $n$ and $s$.  }

\subsection{Preliminary results}

We consider the hypothesis $(H_4)$:

$(H_4)$: there exists a constant $\underline \gamma > 0$ such as $\ds
\underline \gamma \leq \gamma(r)$ for any $r \in (0,1)$.

We obtain the following Lemma (with proof in Annex A):

\begin{Lem}
\label{prelim}
Assume $(H_2), (H_3)$. If $n(.)$ and $s(.)$ are the solutions of the
system $(\cal S_0)$ for a control $e(.)$ then, $ \forall t \in [0,T]$:
(with the convention that the inequalities including $s_0(t)$ are
valid only if $g(r) = g_{\theta}(r)=r^{1-\theta}$)

(i) the function $\ds {g(r(n,s)) \over n}$ satisfies:

$$ {g(r(n^0(t),s^0(t))) \over n^0(t)} \leq {g(r(n(t),s(t))) \over
  n(t)} \leq {g(r(n_0(t),s_0(t))) \over n_0(t)}$$.

Morever, assume $(H_4) $:

(iia) if $\ds 0 < b < {1- {q \over 2} \overline \gamma \over 1- \underline
  \gamma}$ then $\ds n_0(t)s_0(t)^b \leq n(t)s(t)^b \leq n^0(t)s^0(t)^b $

in particular, for $\ds b= {q \over 2}$ the RDI $r(n,s)$ satisfies:

$$ r(n_0(t),s_0(t)) \leq r(n(t),s(t)) \leq r(n^0(t),s^0(t)) $$

(iib) if $\ds b > b_*={1 + {q \over 2} ({1 \over g(r(\underline
    n,s(0)))}-\underline \gamma) \over 1 -\overline \gamma}$ then $\ds
n^0(t)s^0(t)^b \leq n(t)s(t)^b \leq n_0(t)s_0(t)^b $

(iii) the relative increase $\xi$ of the basal area $s$ satisfies $\ds
\xi_m(t) = {s'^0(t) \over s^0(t)^{q \over 2} \overline s^{1-{q \over 2}}}
\leq \xi(t)$. Moreover, if $g(r) = g_{\theta}(r)=r^{1-\theta}$, $\ds
\xi_m(t) = {s'^0(t) \over s^0(t)^{q \over 2} s_0(t)^{1-{q \over 2}}}
\leq \xi(t)$.

\end{Lem}

Remark: If we assume that the final tree-number $n(T)$ is equal to
$\underline n$ then the result obtained in the Lemma \ref{prelim}
remains valid if we replace respectively $n^0,s^0$ by $n_T, s_T$.

\subsection{The optimization problem}

We are focusing here, on the setting in the wood market of a forest
whose evolution is set by the model studied in the previous
paragraphs.  We are introducing the price { (minus the cost of
thinning)} which depends only on the basal area $s$ and the height $h$: we
noted $P_0(s,h,t)$. Owing to the fact that the height  $h$ does not
depend on the basal area $s$ and is equal to a fixed function  $h_0$  of
time $t$, the price can be written in a new function $P$ of $s$ and
$t$: $P(s(t),t):=P_0(s(t),h(t),t)$.

In other words we will set the price function in the following form:
$P_0(s,h,t)= p(s) h e^{-\delta t}$ where $p(.)$ is a price function
for the basal area $s$ and $\delta$ is the actualisation parameter.  We
deduce $P(s,t)= p(s) h_0(t) e^{-\delta t}$ and if we define the
function $\delta_h(.)$ by: $\ds \delta_h(t) = \delta- {h'_0(t) \over
h_0(t)}$ for any $t > 0$ then $P'_t(s,t)=- \delta_h(t) P(s,t)$.
Mostly to simplify we'll assume $p(s) = k s^{\alpha}, \alpha > 0$.

We are assuming that at each time $t$ a quantity $e(t)$ is taken and
that at the end of the period of exploitation $T$ the remaining trees
would have been cut. The instantaneous value of the trees that would
have been cut is $P(s(t), t)e(t)$ and the final value is $P(s(T),
T)n(T)$.

The criterion which we suggest to maximize consists of an integral
term corresponding to the cuttings that would have occurred during the
interval $[0, T]$ and the final term corresponding to the final
cuttings at time $T$.

The optimization problem, relating to the cuttings $e(.)$, on the
interval $[0, T]$, is therefore written:

$$({\mathcal P}): \ \ \ \ \ \ \ \max_{e(.)} \int_0^T P(s(t),t) e(t) dt
 + P(s(T),T) n(T) $$

 with $0 \leq e(t) \leq \overline e$ and $n$ and $s$ solutions of
 $({\mathcal S_0})$ with initial conditions $(n(0), s(0))$ and
 fulfilling the constraints: $n(t) \geq \underline n$ et $r(n(t),s(t))
 \leq 1$.

 Intuitively, from the fact that the function $\ds {g(r(n,s)) \over
   n}$ is a decreasing function of $n$ (Lemma \ref{lemG} (ii)), we are
 tempted to suggest the following assertion:

{\it In order for the trees to get the best benefits from the
nutrients, one should, from the beginning cut a significative number
of trees, so that in the end of the exploitation timescale, one should
get a limited tree number of good quality}.

We will try to validate or invalidate according to the cases this
assertion and we will also try to answer the complementary yet
important questions for management:

1) Does optimal silviculture depend on the term $T$ ?

2) Which role the various parameters of the model play: economic parameters $p(.), \delta$
and growth parameters $g(.)$, $q$ ?

The optimization problem $({\mathcal P})$ can be rewritten just by
replacing $e(t)$ by $- \ds {dn(t) \over dt}$:

$$\max_{n(.) \in {\mathcal C}} -\int_0^T P(s(t),t) {dn(t) \over dt} dt
+P(s(T),T) n(T)$$

 where ${\mathcal C} $ is the whole set of curves:

 $$ {\mathcal C} = \{ n(.) \in C^1([0,T])|-\overline e \leq {dn(t)
   \over dt} \leq 0 \ \& \ A n(t) s(t)^{q \over 2} \leq 1 \}$$

By an integration by part we deduce:

$$\max_{n(.) \in {\mathcal C}} \int_0^T {dP(s(t),t) \over dt} n(t) dt
+P(s(0),0) n(0)$$

  under the same constraints as in the initial problem.


We are here defining the function $\xi(.)$, the relative increase of
the basal area $s$, by $\ds \xi(t) = {s'(t) \over s(t)}$.  
By applying the results of Lemma \ref{prelim} (iii) 
we deduce the following Proposition (with proof in Annex B):

\begin{Pro}
\label{optim}
Assume $(H_1), (H_2), (H_3), (H_4)$, $T \leq \overline T$, then

(i) if $g(r)= r^{1-\theta}$, $\ds \alpha > \alpha_*
=1+(b_*-{q \over 2})(1-\theta)$ and $\ds \delta_h(t) \leq \alpha
(1-\theta) \xi_m(t)$,then the optimal trajectory is $E_0$.

(ii) if $\ds \alpha < {1- {q \over 2} \underline \gamma \over 1
  -\underline \gamma}$ and $\ds \delta_h(t) \leq \alpha \underline
\gamma \xi_m(t)$, then the optimal trajectory is $E^0$.

\end{Pro}

From the remark following the Lemma \ref{prelim} we deduce:

\begin{Cor}
  If we assume that the final tree-number $n(T)$ is equal to
  $\underline n$ then:

  if $\ds \alpha < {1- {q \over 2} \underline \gamma \over 1
    -\underline \gamma}$ and $\ds \delta_h(t) \leq \alpha \underline
  \gamma \xi_m(t)$, then the optimal trajectory is $E_T$.
\end{Cor}

The condition on $\alpha$ in Proposition \ref{optim} (i) implies $p$
must be sufficiently convex. Thus, under the conditions mentioned in
the Proposition \ref{optim} (i) (if $p$ is sufficiently convex and the
parameter of actualization not too high), one may find it beneficial
to cut the maximum tree number at the beginning to ensure a high rate
for the remaining tree basal area at the end of $T$ as foretold in the
stated assertion.  Similar results were obtained with the full model
``Fagacées'' \cite{Loisel}. Moreover, in the studied cases in
Proposition \ref{optim}, silviculture, i.e. cuttings policy $e(.)$,
does not depend on final time $T$. The conditions depend on economic
parameters: a sufficient convexity of the price function relative to
the basal area $s$ and a small enough parameter of actualization.  The
stated assertion however is no longer satisfied under the conditions
of the Proposition \ref{optim} (ii).

\section{Conclusion}

In that article, starting with a tree-growth model governed by the
tree basal area and the number of trees, we study the viability
properties of the system solutions. We highlighted the importance of
the economic parameters and growth parameters on silviculture.  Hence
for a price (minus the thinning costs) is sufficently convex and a
parameter of actualization not too high, it is optimal to cut the
trees at the beginning of the period of exploitation.


\section{Proof of Lemma \ref{prelim}}

(i) the result is a consequence of Lemma \ref{lemG} (ii) (iii) and
Proposition \ref{mini} (ii) (iii).

(iia) $\ds {d(n(t)s^b(t))^a \over dt} = a b n(t)^{a-1} s(t)^{ab-1}
g(r(n(t),s(t))) V(t) \ - a e(t) n(t)^{a-1} s(t)^{ab}$.

If we denote $y(n,s)$ the expression of $\ds {d(n s^b)^a \over dt}$,
then, if we assume $0 < a < 1$:

\begin{align}
  y'_n= & a n(t)^{a-2} s(t)^{ab-1}(b[rg'(r) +(a-1) g(r)](n(t),s(t))
  V(t)  +(1-a)e(t) s(t)) \notag \\
  \geq & a b n(t)^{a-2} s(t)^{ab-1} (\underline \gamma +a-1)
  g(r(n(t),s(t)))  V(t)  \notag \\
  y'_s= & a b n(t)^{a-1} s(t)^{ab-2}([{q \over 2} rg'(r) -(1-ab)
  g(r)](n(t),s(t)) V(t) - a e(t) s(t)) \notag \\
  \leq & a b n(t)^{a-1} s(t)^{ab-2}({q \over 2} \overline \gamma
  +ab-1) g(r(n(t),s(t)))) V(t) \notag
\end{align}

hence, if we choose $a$ such that $\ds 1-\underline \gamma < a <
\min({1- {q \over 2} \overline \gamma \over b},1)$, we deduce $\ds {q
  \over 2} \overline \gamma + a b -1 < 0 < \underline \gamma + a -1$
then $y'_n > 0$ and $y'_s < 0$ and:

$\ds {d (n_0(t) s_0(t)^b)^a \over dt} =y(n_0(t),s_0(t)) \leq {d (n(t)
  s(t)^b)^a \over dt} \leq y(n^0(t),s^0(t))= {d (n^0(t) (s^0(t))^b)^a
  \over dt}$

by integration, we obtain the result. 

(iib) From the expression of $y'_n$ and $y'_s$ and using $e(t) \leq
e_r(s(t),t)$, we deduce that, if $\ds a < 1-\overline \gamma$:

\begin{align}
  y'_n \leq & a n(t)^{a-2} s(t)^{ab-1}(b(\overline \gamma +a-1)
  g(r(n(t),s(t)))  +{q \over 2} (1-a)) V(t) \notag \\
  y'_s \geq & a b n(t)^{a-1} s(t)^{ab-2}(({q \over 2} \underline
  \gamma +ab-1)  g(r(n(t),s(t))) - {q \over 2} a) V(t) \notag 
\end{align}

and, if $\ds b > b_1(a) = {q \over 2} {1-a \over 1-a-\overline \gamma}
{1 \over g(r(\underline n,s(0))} $ then $y'_n < 0$.

if $\ds b > b_2(a) = {q \over 2} {1 \over g(r(\underline n,s(0))} + {1
  - {q \over 2} \underline \gamma \over a}$ then $y'_s > 0$.

To obtain the minimal limit value for $b$, as $b_1$ is increasing in
$a$, and $b_2$ is decreasing in $a$, we choose the value $a$ such that
$b_1(a)=b_2(a)$, this value is $\ds a_*={(1-\overline \gamma)(1-{q
    \over 2} \underline \gamma) \over 1 + {q \over 2} ({\overline
    \gamma \over g(r(\underline n,s(0))}-\underline \gamma)}$ hence
the result if $b > b_*=b_i(a_*)$.

(iii) From Lemma \ref{lemG} (i), $ \ds {g(r) \over r} $ is a decreasing
function of $r$ then:

$\ds {s'(t) \over s(t)} = A {g(r(n(t),s(t))) \over r(n(t),s(t))} {V(t)
  \over s(t)^{1-{q \over 2}}} \geq A {g(r(n^0(t),s^0(t))) \over
  r(n^0(t),s^0(t))} {V(t) \over s_M(t)^{1-{q \over 2}}}$. From $s_M(t)
\leq \overline s$ we deduce the result, if $g(r) =
g_{\theta}(r)=r^{1-\theta}$, $s_M(t) = s_0(t)$.  \hfill $\square$

\section{Proof of Proposition \ref{optim}}

(i) Let's consider the auxiliary problem which consists in maximizing,
at each time $t$, the integrand:

\begin{align} {dP(s(t),t) \over dt} n(t) = & P'_s(s(t),t)
  g(r(n(t),s(t))) V(t) +
  P'_t(s(t),t) n(t) \notag \\
  = & k h_0(t) e^{-\delta t} (\alpha s(t)^{\alpha-1} g(r(n(t),s(t)))
  V(t) - \delta_h(t) n(t) s(t)^{\alpha}) \notag
\end{align}

We denote $y= n s^b$ and $z(y,s,t)$ by: $\ds z(y,s,t) =\alpha
s^{\alpha-1} g(A y s^{{q \over 2}-b}) V(t) - \delta_h(t) y s^{\alpha-b}$ then:

\begin{align}
  z'_y= & \alpha s(t)^{\alpha-1} {[rg'(r)](n(t),s(t)) \over y(t)} V(t)
  - \delta_h(t) s(t)^{\alpha-b}\notag \\
  \geq & s(t)^{\alpha-b-1} (\alpha (1-\theta) {g(r(n(t),s(t)))
    \over n(t)} V(t) - \delta_h(t) s(t)) = s(t)^{\alpha-b} (\alpha
  (1-\theta) \xi(t)- \delta_h(t)) \notag
\end{align}

From $\ds \delta_h(t) < \alpha (1-\theta) \xi_m(t)$ we deduce $z'_y >
0$. Moreover if $b \geq \alpha$:

\begin{align}
  z'_s= & s(t)^{\alpha-2} (\alpha ( [({q \over 2}-b) r
  g'(r)+(\alpha-1) g(r)](n(t),s(t)))V(t)
  -  (\alpha-b)\delta_h(t) y(t) s(t)) \notag \\
  \geq & \alpha s(t)^{\alpha-2} (({q \over 2}-b) (1-\theta) +\alpha-1)
  g(r(n(t),s(t))) V(t) \notag
\end{align}

Due to $\ds \alpha > \alpha_*$ we can choose $b$ such that $\ds
\max(b_*, \alpha) < b < {\alpha-1 \over 1-\theta} +{q \over 2}$. Then
$z'_s > 0$, from Lemma \ref{prelim} (iib) we deduce the result.
 
(ii) We then denote $y= n s^{\alpha}$ and $z(y,s,t)$ by: $\ds z(y,s,t)
=\alpha s^{\alpha-1} g(A y s^{{q \over 2}-\alpha}) V(t) - \delta_h(t)
y$.

\begin{align}
  z'_y= & \alpha s(t)^{\alpha-1} {[rg'(r)](n(t),s(t)) \over y(t)} V(t)
  - \delta_h(t) \notag \\
  \geq & \alpha \underline \gamma {g(r(n(t),s(t))) \over n(t) s(t)}
  V(t) - \delta_h(t) = \alpha \underline \gamma \xi(t)- \delta_h(t)
  \notag
\end{align}

From $\delta_h(t) < \alpha \underline \gamma \xi_m(t)$, $z'_y >
0$. Moreover:

\begin{align}
  z'_s= & \alpha s(t)^{\alpha-2} [ ({q \over 2}-\alpha) r
  g'(r)+(\alpha-1) g(r)] (n(t),s(t)) V(t)\notag \\
  = & \alpha s(t)^{\alpha-2} [( ({q \over 2}-\alpha) \gamma(r)
  +\alpha-1) g(r)](n(t),s(t)) V(t) \notag
\end{align}

From $\ds \alpha < {1- {q \over 2} \underline \gamma \over 1 -
  \underline \gamma}$, we deduce : $\ds z'_s < \alpha s(t)^{\alpha-2}
{ 1 - {q \over 2} \over 1 -\underline \gamma} (\underline \gamma -
\gamma(r)) g(r(n(t),s(t))) V(t) \leq 0$ then $\ds z'_s < 0$. Moreover,
from $z'_y > 0$, Lemma \ref{prelim} (iia) with $b=\alpha$ and $z'_s <
0$ we deduce the result.  \hfill $\square$

\end{document}